\documentclass[12pt]{article}
\usepackage[dvips]{graphicx}
\usepackage{amsmath, amssymb, amsthm}
\def\indep{\mathop{\lower.25em\hbox{$-$}\kern-.9em\shortparallel}\,}

\newcommand{\cS}{{\cal S}}
\newcommand{\cC}{{\cal C}}

\newtheorem{theorem}{Theorem}[section]
\newtheorem{lemma}{Lemma}[section]
\newtheorem{corollary}{Corollary}[section]

\newtheorem{proposition}{Proposition}[section]
\theoremstyle{definition}
\newtheorem{ex}{Example}[section]
\title{Iterative proportional scaling via decomposable submodels for
contingency tables} 

\author{
 Yushi Endo and Akimichi Takemura\\
  Graduate School of Information Science and Technology\\
  University of Tokyo}

\date{November, 2008}

\begin{document}
\maketitle

\begin{abstract}
  We propose iterative proportional scaling (IPS) via decomposable
  submodels for maximizing likelihood function of a hierarchical model
  for contingency tables.  In ordinary IPS the proportional scaling is
  performed by cycling through the members of the generating class of
  a hierarchical model.  We propose to adjust more marginals at each
  step.  This is accomplished by expressing the generating class as a
  union of decomposable submodels and cycling through the decomposable
  models.  We prove convergence of our proposed procedure, if the
  amount of scaling is adjusted properly at each step.  We also
  analyze the proposed algorithms around the maximum likelihood
  estimate (MLE) in detail. Faster
  convergence of our proposed procedure is illustrated by numerical
  examples.
\end{abstract}

\noindent
{\it Keywords and phrases:} \ 
decomposable model, 
hierarchical model, 
$I$-projection,  
iterative proportional fitting, 
Kullback-Leibler divergence.

\section{Introduction}
\label{sec:intro}

Iterative proportional scaling algorithm for contingency tables, 
first proposed by Deming and Stephan  \cite{i0}, 
has been well studied and generalized
by many authors.
Ireland and Kullback \cite{i1} proved
convergence of IPS and Fienberg \cite{i5}
gave a simpler proof of convergence from geometric consideration.
Darroch and Ratcliff \cite{i3} made a
generalization to IPS and its geometrical property was studied
by Csisz\'{a}r \cite{i6}.  Csisz\'{a}r 
\cite{i7}
also gave
a more general proof of convergence 
and justified IPS in a general framework.
Extension of IPS to continuous
case was studied in Kullback \cite{i2} and R\"uschendorf \cite{i4}.
Effective algorithms and implementations of IPS have been also studied
by many authors, including \cite{badsberg}, \cite{i9}, \cite{i11},
\cite{i10}, \cite{i8}. 

In this paper, we propose another generalization of IPS based on
decomposable submodels. Decomposable models or graph decompositions
have been already considered  by Jirou\v{s}ek \cite{i11}, Jirou\v{s}ek and
P\v{r}eu\v{c}il \cite{i10} and Malvestuto \cite{i8}.  However they
used decomposable models for efficient implementation of conventional
IPS in the form of tree-computation.  
Here we use decomposable submodels for generalizing IPS itself.
In our algorithm we adjust a larger set of marginals than the
conventional IPS.  The set of marginals form the generating class of a
decomposable submodel.  By adjusting more marginals, our proposed
algorithm achieves a faster convergence to the maximum likelihood
estimate than the conventional IPS, although at present it seems
difficult to theoretically prove that our procedure is always faster.
We prove convergence of our proposed
procedure, if we adjust the amount of scaling at each step.  We also
analyze in detail the behavior of the proposed algorithms around the
maximum likelihood estimate.  As shown in Section
\ref{sec:experiments} our procedure works well in practice without
adjusting the amount of scaling at each step.

As suggested by a referee, it is an important topic to combine the
idea of the present paper and the tree-computation approach for
efficient implementation of IPS.  Although we do not give a general
result, in Section \ref{sec:experiments} we investigate the
combination in the case of cycle models and show effectiveness of the
combination by numerical experiments.

The organization of this paper is as follows.
In Section \ref{sec:preliminaries} we summarize notations and basic
facts on hierarchical models and decomposable models for multiway
contingency tables. In Section \ref{sec:gen-IPS}  we propose a
generalized IPS via
decomposable submodels, prove its convergence and clarify its behavior
close to the maximum likelihood estimate.
In Section \ref{sec:experiments} we perform some numerical experiments to
illustrate the effectiveness of the proposed procedure.  Some
discussions are given in Section \ref{sec:discussions}.

\section{Preliminaries}
\label{sec:preliminaries}

In this section we summarize notations and preliminary materials on
decomposable models and conventional IPS.

We follow the notation of Lauritzen \cite{st2}.  
Let $\Delta$ denote the set
of variables of a multiway contingency table. 
For each $\delta \in \Delta$, ${\cal
  I}_\delta=\{1,2,\ldots ,I_\delta\}$ denotes the set of levels of
$\delta$.  The set of cells is denoted by ${\cal I}=
\times_{\delta\in\Delta} {\cal I}_\delta$. 
Let $n(i)$ denote the frequency of a cell $i\in {\cal I}$
and let
$n=\sum_{i \in {\cal I}}n(i)$ denote the total sample size.  
Throughout the paper we denote the relative frequency (empirical
distribution) by $r(i)=n(i)/n$.
For a cell $i$ and a subset of variables  $a \subset \Delta$, 
the marginal cell of $i$ for $a$ is denoted by $i_a \in {\cal I}_a =\times
_{\delta \in a} {\cal I}_{\delta}$, %
the marginal of $r$ on $a$ is denoted by $r[a]$,
and the marginal relative frequency of $a$ is denoted by $r(i_a)$.

The {\em generating class} of a hierarchical model is the family of the
variable sets indexing the maximal interaction terms in the
hierarchical model. We denote a hierarchical model with generating
class $\cal C$ by $M({\cal C})$, and call the sets in $\cal C$
the {\em generators} of $M({\cal C})$. A
hierarchical model is a {\em decomposable model} if there exists an ordering
$(C_1, \dots, C_m)$ of its generators that satisfies the running intersection
property:
\begin{center}
(RIP)\quad  For each $j$ $(2\leq j\leq m)$, there exists $k$ $(1\leq k\leq j-1)$,
such that $ C_j\cap (C_1\cup C_2\cup \cdots \cup C_{j-1})\subset C_k$.
\end{center} 
Such an ordering is called a {\it perfect sequence}. 
Given a perfect sequence $(C_1, \dots, C_m)$
of the generators of a decomposable model,  let
\[
S_j = C_j\cap (C_1\cup C_2\cup \cdots \cup C_{j-1}). \qquad (2\leq j\leq m)
\]
If no $S_j$ is the empty set, then the decomposable model is said to be
{\em connected}. If this is the case, then each set $S_j$ is called a {\em separator} 
of the generating class of the decomposable model; moreover, both the
generators of the decomposable model and the separators of its
generating class can be graphically viewed as being the (maximal)
cliques and the minimal vertex separators of a suitable chordal,
connected graph, sometimes called the ``adjacency graph'' of the
generating class of the decomposable model. In what follows, we always
assume that a decomposable model is connected. 
In this paper 
\[
{\cal S}=\{ S_2, \dots, S_m\},
\]
denotes the multiset of separators. The number of times %
a separator $S$ appears in $\cal S$ is called the {\em multiplicity} of $S$.

The MLE of the cell probabilities $\{p(i)\}$ under a hierarchical model
$M({\cal C})$ %
is given by the probability distribution denoted by $p_{\cal C}$ 
that belongs to $M({\cal C})$ and
satisfies the marginality constraints 
\begin{equation}
\label{eq:likelihood-equation}
p_{\cal C}[C] = r[C], \qquad \forall C \in {\cal C}.
\end{equation}
Equivalently, $p_{\cal C}$ 
is the extension of the set of probability
distributions $\{r[C]: C \in {\cal C}\}$ that has the maximum entropy. 
If $M({\cal C})$ is a decomposable model then 
$p_{\cal C}(i)$ has the following product-form
expression:
\begin{equation}
\label{eq:mle}
p_{\cal C}(i)
=\begin{cases}
\displaystyle \frac{\prod _{C\in {\cal C}}r(i_C)}{\prod _{S\in {\cal S}}r(i_S)}, &
\text{if}\ \   r(i_C) >0,\ \forall C\in \cal C,\\
0, & \text{otherwise}.
\end{cases}
\end{equation}
In the following we call $p_{\cal C}$ in (\ref{eq:mle}) 
the {\em maximum-entropy extension} of the set of 
probability distributions $\{r[C]: C \in  {\cal C}\}$.
In Algorithm 2 below, we use the  maximum-entropy extension
of the form (\ref{eq:mle}) of the set 
$\{q[C]: C \in  {\cal C}\}$ even when $q$ is not necessarily
normalized to be a probability distribution.
 
For obtaining MLE for other graphical or hierarchical models we need
some iterative procedure.  The following conventional IPS, cycling
through the elements of the generating class,  is commonly
used for this purpose. In the following let $p^{(t)}(i)$ denote the
estimate of the probability of the cell $i$ at the $t$-th step of
iteration
and let $p^{(t)} = \{ p^{(t)}(i)\}$.

\bigskip
\noindent
{\bf Algorithm 0} \qquad (Conventional IPS)\\
Let $p^{(0)}(i)\equiv 1/\vert {\cal I} \vert$. 
The updating formula is given as 
\begin{equation}
 \label{eq:IPS}
p^{(t+1)}=p^{(t)}\times \frac{r[C]}{p^{(t)}[C]}, 
\end{equation}
where $C=C_j$, \ $j=(t\mod m)+1$.
\bigskip

The Kullback-Leibler divergence (KL-divergence) from a probability 
distribution $p$ to another probability distribution $q$ 
is denoted by 
\[
I(p:q)=\sum _{i\in {\cal I}} p(i) \log \frac{p(i)}{q(i)}. 
\] 
The log sum inequality (Chapter 2 of \cite{cover-thomas}) for
non-negative numbers $a_1,\dots,a_N$ and $b_1, \dots, b_N$ is
\[
\sum_{i=1}^N a_i \log \frac{a_i}{b_i} \ge a \log \frac{a}{b},
\qquad a_i \ge 0, \ b_i \ge 0,  \quad  a=\sum_{i=1}^N a_i, \ b=\sum_{i=1}^N
  b_i,
\label{logsum}
\]
where $a\log \frac {a}{0}=\infty$ if $a>0$, and $0\log 0=0$. 
The equality holds if and only if $a_i/b_i={\rm const}$.

\section{Iterative proportional scaling via decomposable submodels}
\label{sec:gen-IPS}

In this section we propose a generalization of conventional IPS and
study its properties.  At each step of our procedure we update a
larger set of marginals, which form a decomposable submodel.
We prove convergence of our proposed procedure, if the
amount of scaling is adjusted properly at each step.   We also give a
detailed analysis of our procedure when the current estimate is close
to MLE.

\subsection{Proposed algorithms}
\label{subsec:algorithms}

We now describe our proposed procedure.  
A model
$M({\cal C}')$ is a submodel of $M({\cal C})$ if each generator of 
$M({\cal C}')$ is contained in
some generator of $M({\cal C})$. 
Let $\{M({\cal C}_1), \dots, M({\cal C}_u)\}$
be a set of decomposable submodels of 
$M({\cal C})$ such that each generator of $M({\cal C})$
is contained in the generating class  of $M({\cal C}_j)$ for some $j$.
In this case we say that 
$\{M({\cal C}_1), \dots, M({\cal C}_u)\}$ {\em spans} $M({\cal C})$.

In our procedure there is a
problem of normalization as discussed below.  Therefore
we denote the non-normalized estimated cell probability at the $t$-th
step by $q^{(t)}$ and the normalized estimated cell probability by
$p^{(t)}$.

\bigskip\noindent
{\bf Algorithm 1}\qquad  Let $q^{(0)}\equiv 1/\vert {\cal I} \vert$. 
We cycle through $\cC_1, \cC_2, \ldots
,\cC_u$ and for the $t$-th step we update the non-normalized
estimated cell probabilities
as follows 
\begin{equation}
\label{eq:a1}
q^{(t+1)}= q^{(t)} \frac{r_j}{q^{(t)}_j},
\qquad
j= (t\mod u)+1, 
\end{equation}
where  $r_j$ is the maximum-entropy extension of the set of probability
distributions $\{r[C]: C \in  {\cal C}_j\}$ and $q^{(t)}_j$ is 
the maximum-entropy extension of the set of probability
distributions $\{q^{(t)}[C]: C \in {\cal C}_j\}$,  
and the normalized cell probabilities as
\begin{equation}
\label{eq:a1a}
p^{(t+1)}(i)=\frac{q^{(t+1)}(i)}{\sum_{k\in {\cal I}} q^{(t+1)}(k)}.
\end{equation}

\begin{ex}
\label{ex1}
Consider a 4-way contingency table $H\times J\times K\times L$ 
and the following hierarchical model 
with generating class ${\cal C}=\{ \{H,J\}, \{J,K\}, 
\{K,L\},\{H,L\}\}$
 (``4-cycle model'') : 
\[
p_{hjkl}=\exp({a_{hj}+b_{jk}+c_{kl}+d_{hl}}).
\]
By slight abuse of notation write $\Delta=\{H,J,K,L\}$.
 The following ${\cal C}_1$ and ${\cal C}_2$ 
is an example of the family
 of submodels that %
 spans ${\cal C}$.
 \begin{eqnarray*}
  \cC_1=\{ \{ H,J\},\{ J,K\},\{K,L\} \}, \\
  \cC_2=\{ \{ H,J\},\{ K,L\},\{ H,L\} \}. 
 \end{eqnarray*}
For each submodel, the updating procedure is performed as follows. 
\begin{eqnarray*}
q^{(t+1)}(i)
=q^{(t)}(i)\times \frac{r(i_{hj})\times 
r(i_{jk})\times r(i_{kl})}{r(i_j)
\times r(i_k)}
\times \frac{q^{(t)}(i_j)\times q^{(t)}(i_k)}{q^{(t)}(i_{hj})\times q^{(t)}(i_{ jk})\times q^{(t)}(i_{kl})},
\\
q^{(t+2)}(i)
=q^{(t+1)}(i)\times \frac{r(i_{hj})\times r(i_{kl})\times r(i_{hl})}{r(i_h)\times r(i_l)}\times
 \frac{q^{(t+1)}(i_h)\times q^{(t+1)}(i_l)}{q^{(t+1)}(i_{hj})\times q^{(t+1)}(i_{kl})\times q^{(t+1)}(i_{hl})}.
\end{eqnarray*}
\end{ex}

If we set ${\cal C}_1=\{C_1\},\ldots, {\cal C}_m=\{C_m\}$, 
Algorithm 1 coincides with the conventional IPS.
${\cal C}_1,\ldots,{\cal C}_m$ %
span ${\cal C}$. 
Each ${\cal C}_j$ is composed of one generator of the model. %
Hence  $M({\cal C}_j)$ is a decomposable submodel of $M({\cal C})$.
Therefore Algorithm 1 is a %
generalization of conventional IPS.
In the conventional IPS, $p^{(t+1)}(i)$ in (\ref{eq:IPS}) satisfies
$p^{(t+1)}(i_C) = r(i_C)$, which is a likelihood equation
in (\ref{eq:likelihood-equation}).  
But in general $p^{(t+1)}(i)$ in (\ref{eq:a1a}) does not satisfy 
(\ref{eq:likelihood-equation}).  
In other words, 
from a geometric viewpoint of $I$-projection in Csisz\'{a}r 
(\cite{i7},\cite{i6}), the updating rule (\ref{eq:a1}) is not a
projection. We discuss it again in the next section.

In (\ref{eq:a1}), we update $q^{(t)}(i)$.
It should be noted that we have 
$$
q^{(t)}(i)\times \frac{\prod _{C\in
 \cC_j}r(i_C)}{\prod _{S\in \cS_j}r(i_S)}\times \frac{\prod _{S\in
 \cS_j}q^{(t)}(i_S)}{\prod _{C\in \cC_j}q^{(t)}(i_C)}
 =p^{(t)}(i) \times \frac{\prod _{C\in \cC_j}r(i_C)}{\prod _{S\in
 \cS_j}r(i_S)}\times \frac{\prod _{S\in 
 \cS_j}p^{(t)}(i_S)}{\prod _{C\in \cC_j}p^{(t)}(i_C)}\notag\\ 
$$
because the normalizing constant is canceled on the right-hand side of
(\ref{eq:a1}).   
Also it is easy to see that, if 
Algorithm 1 in terms of  $\{q^{(t)}(i)\}$ converges,
then the limiting $q$'s are automatically normalized.

Unfortunately it is difficult to prove convergence of Algorithm
1. %
The difficulty lies in the fact that the sum $\sum
_{i \in{\cal I}} q^{(t+1)}(i)$ after updating might exceed 1 
(i.e.\ $\sum _{i\in {\cal I}} q^{(t+1)}(i)>1$)
in Algorithm 1 even if $q^{(t)}$ is normalized as 
$\sum _i q^{(t)}(i)=1$. 
However we recommend it because in practice, it works well and has
converged to MLE in all of our experiments and converges faster than the
conventional IPS as shown in Section \ref{sec:experiments}. 

In order to deal with the theoretical difficulty concerning
the normalization of $q^{(t+1)}$ we consider
adjusting the amount of updating. 
At this point, we need the following lemma. 

\begin{lemma}
\label{lem1}
Let $r$ and $q$ be two probability distributions over $\cal I$, and
$M({\cal C}')$ a decomposable model over a nonempty (proper or
improper) subset $\Delta'$ of $\Delta Δ$.  Let $r_{{\cal C}'}$ be the
maximum-entropy extension of the set of probability distributions
$\{r[C]: C \in {\cal C}'\}$ and $q_{{\cal C}'}$ be the maximum-entropy
extension of the set of probability distributions $\{q[C]: C \in {\cal
  C}'\}$. If $q[\Delta ']$ is not an extension of the set of
probability distributions $\{r[C]: C \in {\cal C}'\}$, then there
exists a unique $\alpha_0 \ge 0$ for which the function 
\[
q'=q \times \left(\frac{r_{{\cal C}'}}{q_{{\cal C}'}}\right)^\alpha
\]
is a probability distribution.
\end{lemma}

\begin{proof}
In view of (\ref{eq:mle}) we have
\[
1 = \sum_{i\in {\cal I}} r_{{\cal C}'}(i)=\sum_{i\in {\cal I}} q_{{\cal C}'}(i)
\]
Therefore if 
\begin{equation}
\label{eq:ratio1}
\frac{r_{{\cal C}'}(i)}{q_{{\cal C}'}(i)}\le 1
\end{equation}
for all $i$, then  the equality in 
(\ref{eq:ratio1}) holds for all $i$ with  $q(i) > 0$.
Therefore under the condition of the lemma there exists at least one
cell $i\in {\cal I}$ such that
\[
\frac{r_{{\cal C}'}(i)}{q_{{\cal C}'}(i)} > 1, \qquad q(i)>0.
\]
Then $q(i) (r_{{\cal C}'(i)}/q_{{\cal C}'}(i))^\alpha$ %
for this $i$ is
strictly convex in $\alpha$ and diverges to $+\infty$ as $\alpha
\rightarrow \infty$.  
Write 
\begin{equation}
\label{eq:g-alpha}
g(\alpha) = 
\sum_{i\in {\cal I}} q(i)
 \left(\frac{r_{{\cal C}'}(i)}{q_{{\cal C}'}(i)}\right)^\alpha.
\end{equation}
Then %
$g(\alpha)$%
is also strictly convex in $\alpha$ and diverges to $+\infty$ as
$\alpha \rightarrow \infty$.

Write ${\cal C}' = \{ C_1,\dots, C_v \}$ and 
${\cal S}'= \{ S_2,\dots, S_v \}$.
Consider the differential of $g(\alpha)$ at $\alpha=0$.
\begin{align*}
  g'(0)&=\sum_{i} q(i) \log \left( \frac{\prod _{C\in
        \cC'}r(i_C)}{\prod _{S\in \cS'}r(i_S)}
    \times \frac{\prod _{S\in \cS'}q(i_S)}
           {\prod _{C\in \cC'}q(i_C)}\right)\\
  &=\sum _{C\in \cC'} \sum_i q(i) \log
  \frac{r(i_C)}{q(i_C)}
  - \sum _{S\in \cS'} \sum_i q(i) \log\frac{r(i_S)}{q(i_S)}\\
  &=\sum_{C\in \cC'} \sum_{i_C} q(i_C) \log
  \frac{r(i_C)}{q(i_C)} - \sum_{S\in \cS'} \sum_{i_S}
  q(i_S) \log\frac{r(i_S)}{q(i_S)}.
\end{align*}
Now %
\[
\sum_{i_{C_1}} q(i_{C_1}) 
\log \frac{r(i_{C_1})}{q(i_{C_1})}
\]
is the negative of KL-divergence  and nonpositive. By the log sum inequality, 
\[
\sum_{i_{C_k}} q(i_{C_k}) \log \frac{r(i_{C_k})}{q(i_{C_k})}
- \sum_{i_{S_k}} q^{(t)}(i_{S_k}) \log \frac{r(i_{S_k})}{q(i_{S_k})}
\]
is also nonpositive for $2\leq k\leq v$. Equality holds if and only if 
\[
 r(i_{C})=q^{(t) }(i_{C}), \quad \forall C\in \cC', \ \forall i_C.
\]
Then, except for such a case, $g(0)=1$, $g'(0) < 0$,
$g(\infty)=\infty$, and $g(\alpha)$ is strictly convex in
 $\alpha$.  
Therefore there exists a unique $\alpha_0 > 0$ 
such that $g(\alpha_0) =1$.
\end{proof}

We now present the following algorithm  and its modification
based on Lemma  \ref{lem1}.

\bigskip
\noindent {\bf Algorithm 2}\qquad  Let %
$\alpha^{(t)}\ge 0$.
We cycle through $\cC_1, \cC_2, \ldots ,\cC_u$ and for the $t$-th step
we update the unnormalized estimated cell probabilities as
\begin{equation}
  q^{(t+1)}
  =q^{(t)}\times \left( \frac{r_j}{q^{(t)}_j}
\right)^{\alpha^{(t)}}, \qquad
j= (t\mod u)+1, 
\label{eq:algorithm2}
\end{equation}
and the normalized cell probabilities as
$p^{(t+1)}(i)=q^{(t+1)}(i)/\sum_{k\in {\cal I}} q^{(t+1)}(k)$.

\bigskip
Note that also in Algorithm 2 we do not need to normalize at each step and 
we can perform normalization any time, because $\{ q^{(t)}(i)\}$ is
always proportional to $\{ p^{(t)}(i)\}$.

\bigskip

\noindent {\bf Algorithm 3}\qquad  We cycle through $\cC_1, \cC_2, \ldots
,\cC_u$ and for the $t$-th step we update the estimated cell probabilities
as follows 
\begin{equation}
\label{eq:a3}
  q^{(t+1)}
  =q^{(t)}\times \left( \frac{r_j}{q^{(t)}_j}
\right)^{\alpha^{(t)}_0}, \qquad
j= (t\mod u)+1, 
\end{equation}
where $\alpha^{(t)}_0 \geq 0$ is given in Lemma \ref{lem1} with
$q=q^{(t)}$.

\subsection{Correctness of the proposed algorithms}
\label{subsec:convergence}

In this section, we prove the correctness %
of proposed algorithms. 
As before let 
$\{r(i)\}$ denote the empirical distribution and let $\{p_{\cal C}(i)\}$ 
denote the MLE.  Because we consider hierarchical models, the following
equation holds (\cite{i7}, \cite{i6}).
\[
I(r:q)=I(r:p_{\cal C})+I(p_{\cal C}:q).
\]
$I(r:q)$ corresponds to the log likelihood.  Therefore we can prove
the %
correctness of our algorithms by proving 
$I(p_{\cal C}:q^{(t)}) \rightarrow 0$ as $t \rightarrow\infty$.

\begin{theorem} 
{\rm Algorithm 3} converges to MLE.
\label{t1}
\end{theorem}

\begin{proof} 
Consider KL-divergence after updating, 
\begin{eqnarray*}
  I(p_{\cal C};p^{(t+1)})
&=&\sum _{i}p_{\cal C}(i)\log \frac{p_{\cal C}(i)}{p^{(t+1)}(i)} \\
&=&\sum _{i}p_{\cal C}(i)\log \frac{p_{\cal C}(i)}{p^{(t)}(i)\times
 \left(  \frac{\prod _{C\in \cC_j}r(i_C)}{\prod _{S\in \cS_j}r(i_S)}
  \times \frac{\prod _{S\in \cS_j}p^{(t)}(i_S)}{\prod _{C\in
      \cC_j}p^{(t)}(i_C)}\right) ^{\alpha^{(t)} _0}}
\\
&=&\sum _{i}p_{\cal C}(i)\log \frac{p_{\cal C}(i)}{p^{(t)}(i)}
 - \alpha^{(t)} _0 \sum _{i}p_{\cal C}(i)\log \frac{\prod _{C\in
     \cC_j}r(i_C)}{\prod _{S\in \cS_j}r(i_S)}
   \times \frac{\prod_{S\in \cS_j}p^{(t)}(i_S)}
  {\prod_{C\in \cC_j}p^{(t)}(i_C)}. 
\end{eqnarray*}
Write ${\cal C}_j = \{ C_1,\dots, C_v \}$ and 
${\cal S}_j = \{ S_2,\dots, S_v \}$ as in the proof of Lemma
\ref{lem1}.
Then, 
\[
\sum_{i} p_{\cal C}(i) \log
\frac{r(i_{C_1})}{p^{(t)}(i_{C_1})}=\sum_{i_{C_1}} r(i_{C_1}) \log
\frac{r(i_{C_1})}{p^{(t)}(i_{C_1})}
\]
is a KL-divergence, and nonnegative. By the log sum inequality, 
\[
\sum_{i_{C_k}} r(i_{C_k})
\log \frac{r(i_{C_k})}{p^{(t)}(i_{C_k})}-\sum_{i_{S_k}} r(i_{S_k}) \log
\frac{r(i_{S_k})}{p^{(t)}(i_{S_k})}
\]
is also nonnegative for $2\leq k\leq v$. Therefore, 
\[
I(p_{\cal C};p^{(t+1)})\leq I(p_{\cal C};p^{(t)})
\]
holds. Equality holds if and only if 
$r(i_C)=q^{(t) }(i_C)$, $\forall  C\in \cC_j$.
We see that $I(p_{\cal C};p^{(t)})$ always decreases after updating. The
rest of the proof is the same as the classical one (\cite{i2}).
\end{proof}

\begin{corollary}
  Using $0< \alpha^{(t)} \leq \alpha_0^{(t)}$, {\rm Algorithm 2}
  converges to MLE.
\end{corollary}
\begin{proof}
Consider KL-divergence after updating, 
\begin{align*}
& I(p_{\cal C};p^{(t+1)})\\
& \quad = \sum _{i}p_{\cal C}(i)\log \frac{p_{\cal C}(i)}{q^{(t)}(i)}
 - \alpha^{(t)} \sum_{i}p_{\cal C}(i)\log \frac{\prod_{C\in
     \cC_j}r(i_C)}{\prod_{S\in \cS_j}r(i_S)}
   \times \frac{\prod_{S\in \cS_j}q^{(t)}(i_S)}{\prod_{C\in \cC_j}q^{(t)}(i_C)}
+ \log g(\alpha^{(t)}),
\end{align*}
where $g(\alpha)$ is given in (\ref{eq:g-alpha}) with $q=q^{(t)}$.
Because $\alpha^{(t)} \leq \alpha^{(t)}_0$, $\log g(\alpha^{(t)})$ 
is nonpositive and $I(p_{\cal C};p^{(t)})$ always decreases after
updating. The rest of the proof is the same as Theorem \ref{t1}. 
\end{proof}

At this point we discuss Algorithm 3 from a geometric viewpoint of
$I$-projection in the sense of Csisz\'{a}r (\cite{i7}, \cite{i6}).  In
our procedure we adjust a larger set of marginals than the
conventional IPS and in practice KL-divergence decreases more in
our proposed algorithms than the conventional IPS for each step.  
However it is difficult to guarantee this theoretically.  
The difficulty lies in the fact that 
the updating rule (\ref{eq:a3}) is not a projection.
In fact,
if we repeat (\ref{eq:a3}) twice with the same $\cC_j$ then the cell
probabilities change, whereas in the conventional IPS repeating the
same updating step twice does not change the cell probabilities after
the first update.  We can understand the situation as follows.  
Starting from the current
estimate $\{ p^{(t)}(i)\}$ suppose 
that we repeat the step (\ref{eq:a3}) with the same $\cC_j$ until the cell
probabilities converge to $\{ p^{\star}(i)\}$.  
Then the limit $\{ p^{\star}(i)\}$ maximizes the likelihood function
among $\{p(i)\}$ of the form
\begin{equation}
\label{eq:log-affine}
p(i) = p^{(t)}(i) \prod_{C\in \cC_j} \mu(i_C).
\end{equation}
The right-hand side of (\ref{eq:log-affine}) forms a 
log-affine model through  $\{p^{(t)}(i)\}$ (Section 4.2.3 of \cite{st2}).
Since updating a single $C \in \cC_j$ in the conventional IPS
is a special case of (\ref{eq:log-affine}), it follows that 
\begin{equation}
\label{eq:pstar}
I(p_{\cal C}:p^{\star}) \le I(p_{\cal C}:p^{(t+1)'}),
\end{equation}
where 
$\{p^{(t+1)'}(i)\}$ is the updated estimate by the conventional IPS
for some $C\in \cC_j$. Therefore a larger decrease of KL-divergence
of our procedure compared to conventional IPS is only guaranteed 
in the sense of (\ref{eq:pstar}).  The situation will become more
clear when we analyze the behavior of Algorithm 3 close to MLE in the
next section.

\subsection{Analysis of behavior close to the maximum likelihood estimate}
\label{subsec:local-behavior}
In this section, we study the behavior of our algorithms when the 
current estimate is already close to MLE.  
We assume that MLE is in the
interior of the parameter space and $p_{\cal C}(i)>0$ for all
$i\in {\cal I}$. 
We analyze the behavior of $\alpha^{(t)}_0$. We also consider the value of $\alpha^{(t)}=\alpha^{(t)}_1$ which
reduces the  KL-divergence most and the value of $\alpha^{(t)}=\alpha^{(t)}_2$
such that 
KL-divergence decreases in Algorithm 2 for $0\le \alpha^{(t)} \le \alpha^{(t)}_2$.

We repeatedly use the following expansion, 
\begin{equation}
\log (1+x) = x -\frac{x^2}{2}+O(x^3), \quad x\rightarrow 0.
\end{equation}
Assume that the current estimate $\{p^{(t)}(i)\}$
is close to MLE in the following sense.
For sufficiently small $\varepsilon >0$ and for all $C\in \cC$, $S\in
\cS$, $i_C$, $i_S$  we have
\begin{equation}
1-\varepsilon < \frac{r(i_C)}{p^{(t)}(i_C)},\frac{r(i_S)}{p^{(t)}(i_S)} <1+\varepsilon .
\label{ep}
\end{equation}

The following proposition describes the behavior of $\alpha^{(t)}_0$ in
Algorithm 3.

\begin{proposition}
\label{prop:3.1}
Assume 
$\{p^{(t)}(i)\}$
is close to MLE  in the sense of (\ref{ep}).
Then 
\begin{equation}
\alpha^{(t)}_0
=\frac{\sum _i p^{(t)}(i)\left\{\sum _{C\in \cC_j}
\left( \frac {r(i_C)}{p^{(t)}(i_C)}-1\right)^2 
- \sum _{S\in \cS_j} \left( \frac {r(i_S)}{p^{(t)}(i_S)}-1\right)^2
\right\}}
{\sum_{i} p^{(t)}(i) 
\left\{ \sum _{C\in \cC_j}\left( \frac {r(i_C)}{p^{(t)}(i_C)}-1
  \right) 
- \sum _{S\in \cS_j}\left( \frac {r(i_S)}{p^{(t)}(i_S)}-1\right)
\right\}^2 }
+ O(\varepsilon ).
\label{eq:alpha0-close}
\end{equation}
\label{alpha0}
\end{proposition}

Before giving a proof of this Proposition we rewrite the numerator of
the right-hand side of (\ref{eq:alpha0-close}).
Let ${\cal C}_j = \{ C_1,\dots, C_v \}$ and 
${\cal S}_j = \{ S_2,\dots, S_v \}$. %
Then
\begin{align}
&\sum _i p^{(t)}(i)\left\{\sum _{C\in \cC_j}
\left( \frac {r(i_C)}{p^{(t)}(i_C)}-1\right)^2 
- \sum _{S\in \cS_j} \left( \frac {r(i_S)}{p^{(t)}(i_S)}-1\right)^2
\right\}
\nonumber \\
& \qquad
= \sum_{i_{C_1}} p^{(t)}(i_{C_1}) \left(
\frac {r(i_{C_1})}{p^{(t)}(i_{C_1})}-1\right)^2 
\nonumber \\
& \qquad \qquad
+ \sum_{k=2}^v \sum_{i_{C_k}}
p^{(t)}(i_{C_k}) \left( \frac{r(i_{C_k})}{p^{(t)}(i_{C_k})} - 
\frac{r(i_{S_k})}{p^{(t)}(i_{S_k})}\right)^2.
\label{eq:prop31-numerator}
\end{align}
Therefore the numerator is nonnegative.  Also note that the
denominator of the right-hand side of 
(\ref{eq:alpha0-close}) can be written as
\begin{align}
&\sum _{C\in \cC_j}\left( \frac {r(i_C)}{p^{(t)}(i_C)}-1
  \right) 
- \sum _{S\in \cS_j}\left( \frac {r(i_S)}{p^{(t)}(i_S)}-1\right)
\nonumber \\ & \qquad \qquad
=\left(
\frac {r(i_{C_1})}{p^{(t)}(i_{C_1})}-1\right)
+ \sum_{k=2}^v 
\left( \frac{r(i_{C_k})}{p^{(t)}(i_{C_k})} - 
\frac{r(i_{S_k})}{p^{(t)}(i_{S_k})}\right).
\label{eq:prop31-denom}
\end{align}
We see that the numerator of $\alpha^{(t)}_0$ consists of the diagonal square terms
when we expand the square of denominator in the form of 
(\ref{eq:prop31-denom}).  
We now give a proof of Proposition
\ref{prop:3.1}.

\begin{proof}
Consider the following expansion, 
\begin{align*}
\log \frac{\prod _{C\in \cC_j}r(i_C)}{\prod _{S\in \cS_j}r(i_S)}
\times \frac{\prod _{S\in \cS_j}p^{(t)}(i_S)}{\prod _{C\in \cC_j}p^{(t)}(i_C)}
=&\sum _{C\in \cC_j}\log \frac {r(i_C)}{p^{(t)}(i_C)}-\sum _{S\in \cS_j}\log \frac {r(i_S)}{p^{(t)}(i_S)}\displaybreak[0] \\
=&\sum _{C\in \cC_j}\left( \frac {r(i_{C_j})}{p^{(t)}(i_C)}-1 \right)-\sum _{S\in \cS_j}\left( \frac {r(i_S)}{p^{(t)}(i_S)}-1\right)+O(\varepsilon ^2)\displaybreak[0] \\
=&O(\varepsilon ). 
\end{align*}
Then the $s$-th derivative of $g(\alpha^{(t)} )$ at $0$ is
\begin{align*}
g^{(s)}(0)=&\sum_{i} p^{(t)}(i) \left( \log \frac{\prod _{C\in \cC_j}r(i_C)}{\prod _{S\in \cS_j}r(i_S)}
\times \frac{\prod _{S\in \cS_j}p^{(t)}(i_S)}{\prod _{C\in \cC_j}p^{(t)}(i_C)} \right)^s\\
=&O(\varepsilon ^s).
\end{align*}
The first and the second order derivatives of $g(\alpha^{(t)} )$ at $0$ are, 
\begin{align*}
g^{(1)}(0)=&\sum_{i} p^{(t)}(i) \left( \log \frac{\prod _{C\in \cC_j}r(i_C)}{\prod _{S\in \cS_j}r(i_S)}
\times \frac{\prod _{S\in \cS_j}p^{(t)}(i_S)}{\prod _{C\in \cC_j}p^{(t)}(i_C)} \right)\\
=&\sum _{C\in \cC_j}\sum _{i_C} p^{(t)}(i_C)\log \frac {r(i_C)}{p^{(t)}(i_C)}-\sum _{S\in \cS_j}\sum _{i_S} p^{(t)}(i_S)\log \frac {r(i_S)}{p^{(t)}(i_S)} \displaybreak[0] \\
=&\sum _{C\in \cC_j}\sum _{i_C} p^{(t)}(i_C)\left\{ \frac {r(i_C)}{p^{(t)}(i_C)}-1-\frac{1}{2}\left( \frac {r(i_C)}{p^{(t)}(i_C)}-1\right)^2\right\}\\
&-\sum _{S\in \cS_j}\sum _{i_S} p^{(t)}(i_S)\left\{ \frac {r(i_S)}{p^{(t)}(i_S)}-1-\frac{1}{2}\left( \frac {r(i_S)}{p^{(t)}(i_S)}-1\right)^2 \right\}+O(\varepsilon ^3)\displaybreak[0] \\
=&\sum _{C\in \cC_j}\sum _{i_C} \left\{  (r(i_C)-p^{(t)}(i_C))- \frac{p^{(t)}(i_C)}{2}\left( \frac {r(i_C)}{p^{(t)}(i_C)}-1\right)^2\right\}\\
&-\sum _{S\in \cS_j}\sum _{i_S} \left\{ (r(i_S)-p^{(t)}(i_S))- \frac{p^{(t)}(i_S)}{2}\left( \frac {r(i_S)}{p^{(t)}(i_S)}-1\right)^2 \right\}+O(\varepsilon ^3)\displaybreak[0] \\
=&\frac{1}{2}\sum _i p^{(t)}(i)\left\{ \sum _{S\in \cS_j} \left( \frac {r(i_S)}{p^{(t)}(i_S)}-1\right)^2 -\sum _{C\in \cC_j}\left( \frac {r(i_C)}{p^{(t)}(i_C)}-1\right)^2\right\} +O(\varepsilon ^3),
\end{align*}
and 
\begin{align*}
g^{(2)}(0)=&\sum_{i} p^{(t)}(i) \left( \log \frac{\prod _{C\in \cC_j}r(i_C)}{\prod _{S\in \cS_j}r(i_S)}
\times \frac{\prod _{S\in \cS_j}p^{(t)}(i_S)}{\prod _{C\in \cC_j}p^{(t)}(i_C)} \right)^2 \displaybreak[0]\\
=&\sum_{i} p^{(t)}(i) \left\{ \sum _{C\in \cC_j}\left( \frac{r(i_C)}{p^{(t)}(i_C)}-1 \right) -\sum _{S\in \cS_j}\left( \frac{r(i_S)}{p^{(t)}(i_S)}-1\right) \right\}^2 
+O(\varepsilon^3).
\end{align*}
Then, we expand $g(\alpha^{(t)} )$ at $0$, 
\[
g(\alpha^{(t)} ) = g(0)+\alpha^{(t)}
 g^{(1)}(0)+\frac{(\alpha^{(t)})^2}{2}g^{(2)}(0)+O(\varepsilon  ^3). 
\]
Assuming normalization at each step of the algorithm, we have $g(0)=1$ and substituting 
$\alpha^{(t)} _0$ for  $\alpha^{(t)}$, we obtain
\begin{align*}
\alpha^{(t)} _0
&=\frac{-2g^{(1)}(0)}{g^{(2)}(0)}+ O(\varepsilon) \\
&=\frac{\sum _i p^{(t)}(i)\left\{\sum _{C\in \cC_j}\left( \frac
      {r(i_C)}{p^{(t)}(i_C)}-1\right)^2 
   - \sum _{S\in \cS_j} \left( \frac {r(i_S)}{p^{(t)}(i_S)}-1\right)^2
 \right\}}{\sum_{i} p^{(t)}(i) \left\{ \sum _{C\in \cC_j}\left( \frac
     {r(i_C)}{p^{(t)}(i_C)}-1 \right) 
  -\sum _{S\in \cS_j}\left( \frac {r(i_S)}{p^{(t)}(i_S)}-1\right) 
   \right\}^2 }+O(\varepsilon ).
\end{align*}
\end{proof}

Consider (\ref{eq:prop31-numerator})
and (\ref{eq:prop31-denom}).  If the signs of  the terms on the right
hand side of (\ref{eq:prop31-denom}) are ``random'' 
then we can expect that $\alpha^{(t)}_0$ is close to 1.
We can imagine that $\{p^{(t)}(i)\}$ converges to MLE from various
directions. Then $\alpha^{(t)}_0$ is close to 1 ``on the average''.
Furthermore as shown in the following proposition $\alpha^{(t)}_0$ is the
optimum value of the adjustment close to MLE. 
We believe that this is the reason that Algorithm 1 works
very well in practice.

\begin{proposition}
Assume 
$\{p^{(t)}(i)\}$
is close to MLE  in the sense of (\ref{ep}). Then
\begin{equation}
\alpha^{(t)} _1=\alpha^{(t)} _0+O(\varepsilon ), 
\end{equation}
where $\alpha^{(t)} _1$ is the value of $\alpha^{(t)}$ 
which reduces the KL-divergence most.
\label{alpha1}
\end{proposition}
\begin{proof}
Define $F(\alpha^{(t)})$ by 
\begin{equation}
\label{eq:F-alpha}
F(\alpha^{(t)})=\alpha^{(t)} \sum_{i} p^{*}(i)  \log \frac{\prod _{C\in \cC_j}r(i_C)}{\prod _{S\in \cS_j}r(i_S)}\times \frac{\prod _{S\in \cS_j}p^{(t)}(i_S)}{\prod _{C\in \cC_j}p^{(t)}(i_C)},
\end{equation}
which corresponds to the decrease of KL-divergence before
normalization.
Consider the derivative of $F(\alpha^{(t)})$, 
\begin{align*}
F^{(1)}(\alpha^{(t)})=&\sum_{i} p^{*}(i)  \log \frac{\prod _{C\in \cC_j}r(i_C)}{\prod _{S\in \cS_j}r(i_S)}\times \frac{\prod _{S\in \cS_j}p^{(t)}(i_S)}{\prod _{C\in \cC_j}p^{(t)}(i_C)} \displaybreak[0] \\
=&g^{(1)}(0)+ \sum_{i} p^{(t)}(i) \left(\frac{p^{*}(i)}{p^{(t)}(i)}-1 \right) \log \frac{\prod _{C\in \cC_j}r(i_C)}{\prod _{S\in \cS_j}r(i_S)}\times \frac{\prod _{S\in \cS_j}p^{(t)}(i_S)}{\prod _{C\in \cC_j}p^{(t)}(i_C)}\displaybreak[0] \\
=&g^{(1)}(0)+ \sum _i p^{(t)}(i)\left(\frac{p^{*}(i)}{p^{(t)}(i)}-1 \right)\left\{\sum _{C\in \cC_j}\left( \frac {r(i_C)}{p^{(t)}(i_C)}-1 \right)- \sum _{S\in \cS_j} \left( \frac {r(i_S)}{p^{(t)}(i_S)}-1\right) \right\}\\
& +O(\varepsilon ^3) \displaybreak[0] \\
=&g^{(1)}(0)+ \sum _i p^{(t)}(i)\left\{ \sum _{C\in \cC_j}\left( \frac {r(i_C)}{p^{(t)}(i_C)}-1\right)^2-\sum _{S\in \cS_j} \left( \frac {r(i_S)}{p^{(t)}(i_S)}-1\right)^2 \right\} +O(\varepsilon ^3) \displaybreak[0] \\
=&-g^{(1)}(0)+O(\varepsilon ^3).
\end{align*}
Consider the derivative of $F(\alpha^{(t)})-\log g(\alpha^{(t)})$ and equating 0, we obtain, 
\[
F^{(1)}(\alpha^{(t)} _1)-\frac{g^{(1)}(\alpha^{(t)} _1)}{g(\alpha^{(t)} _1)}=0.
\]
Then
\begin{align*}
g^{(1)}(\alpha^{(t)} _1)&=g^{(1)}(0)+\alpha^{(t)} _1 g^{(2)}(0)+O(\varepsilon ^3),\\
g(\alpha^{(t)} _1)&=g(0)+\alpha^{(t)} _1
 g^{(1)}(0)+\frac{(\alpha^{(t)}_1)^2}{2}g^{(2)}(0)+O(\varepsilon ^3)\\ 
&=1+O(\varepsilon ^2)
\end{align*}
and 
\[
F^{(1)}(\alpha^{(t)} _1)-g^{(1)}(\alpha^{(t)} _1)+O(\varepsilon ^3)=-g^{(1)}(0)-g^{(1)}(0)-\alpha^{(t)} _1 g^{(2)}(0)+O(\varepsilon ^3)=0.
\]
Therefore we have
\[
\alpha^{(t)}_1=\frac{-2g^{(1)}(0)}{g^{(2)}(0)}+O(\varepsilon )
=\alpha^{(t)} _0+O(\varepsilon ).
\]
\end{proof}

Finally we show that KL-divergence decreases 
in the range $0 < \alpha^{(t)} < 2\alpha^{(t)}_0$.  This result indicates that 
in Algorithm 2, $\alpha^{(t)} > \alpha^{(t)}_0$ often decreases KL-divergence in
practice.

\begin{figure}[htb]
\begin{center}
\begin{minipage}{110mm}{
\begin{center}
    \includegraphics[scale=0.65]{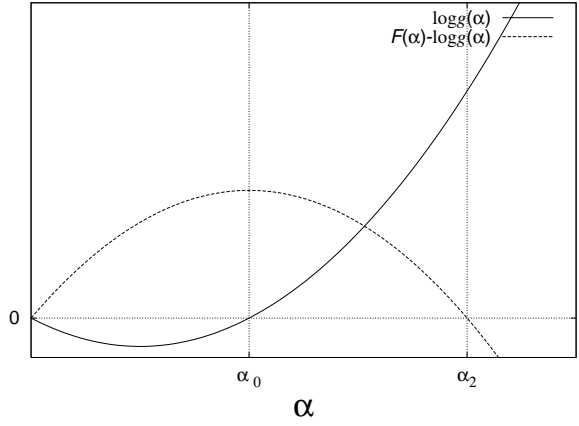}
   \caption{Behavior of $\log g(\alpha^{(t)})$ and $F(\alpha^{(t)}
 )-\log g(\alpha^{(t)})$}
\label{behavior}
\end{center}}
\end{minipage}
\end{center}
\end{figure}

\begin{proposition}
Assume 
$\{p^{(t)}(i)\}$
is close to MLE  in the sense of (\ref{ep}). Then
\begin{equation}
\alpha^{(t)} _2=2\alpha^{(t)} _0+O(\varepsilon ).
\end{equation}
where $\alpha^{(t)} _2$ is the value of $\alpha^{(t)}$ such that 
$I(p_{\cal C}: p^{(t+1)}) = I(p_{\cal C}: p^{(t)})$ in {\rm Algorithm 2}.
\label{alpha2}
\end{proposition}
\begin{proof}
\[
0=F(\alpha^{(t)} _2)-\log g(\alpha^{(t)} _2)
=-\alpha^{(t)} _2g^{(1)}(0) -\alpha^{(t)}
_2g^{(1)}(0)+\frac{(\alpha^{(t)}_2)^2}{2}g^{(2)}(0)+O(\varepsilon ^3)
\]
and
\[
\alpha^{(t)} _2=\frac{-4g^{(1)}(0)}{g^{(2)}(0)}+O(\varepsilon)
=2\alpha^{(t)} _0+O(\varepsilon).
\]
\end{proof}
We show the behavior of $\log g(\alpha^{(t)})$ and $F(\alpha^{(t)} )-\log g(\alpha^{(t)})$ in
Figure \ref{behavior}. 
Proposition \ref{alpha0}, Proposition \ref{alpha1} and
Proposition \ref{alpha2} indicate that in many cases we can decrease KL-divergence
by using $\alpha^{(t)} =1$. In the next section we illustrate this by numerical experiments. 

\section{Numerical experiments for cycle models}
\label{sec:experiments}
\begin{figure}[th]
\setlength{\unitlength}{1.0mm}
\begin{minipage}{40mm}{
\begin{center}
\begin{picture}(33,26)
\put(11,14){\circle{5}}
\put(11,5){\circle{5}}
\put(29,5){\circle{5}}
\put(29,14){\circle{5}}
\put(20,23){\circle{5}}
\put(10,12.8){1}
\put(10,3.8){2}
\put(28,3.8){3}
\put(28,12.8){4}
\put(19,21.8){5}
\put(11,7.5){\line(0,1){4}}
\put(13.5,5){\line(1,0){13}}
\put(29,7.5){\line(0,1){4}}
\put(12.77,15.77){\line(1,1){5.46}}
\put(27.23,15.27){\line(-1,1){5.46}}
\end{picture}
\end{center}}
\end{minipage}
\begin{minipage}{5mm}{
\begin{picture}(10,10)
\put(5,5){=}
\end{picture}
}
\end{minipage}
\begin{minipage}{40mm}{
\begin{center}
\begin{picture}(33,26)
\put(11,14){\circle{5}}
\put(11,5){\circle{5}}
\put(29,5){\circle{5}}
\put(29,14){\circle{5}}
\put(20,23){\circle{5}}
\put(10,12.8){1}
\put(10,3.8){2}
\put(28,3.8){3}
\put(28,12.8){4}
\put(19,21.8){5}
\put(11,7.5){\line(0,1){4}}
\put(13.5,5){\line(1,0){13}}
\put(29,7.5){\line(0,1){4}}
\put(27.23,15.27){\line(-1,1){5.46}}
\end{picture}
\end{center}}
\end{minipage}
\begin{minipage}{5mm}{
\begin{picture}(10,10)
\put(5,5){$\cup$}
\end{picture}
}
\end{minipage}
\begin{minipage}{40mm}{
\begin{center}
\begin{picture}(33,26)
\put(11,14){\circle{5}}
\put(11,5){\circle{5}}
\put(29,5){\circle{5}}
\put(29,14){\circle{5}}
\put(20,23){\circle{5}}
\put(10,12.8){1}
\put(10,3.8){2}
\put(28,3.8){3}
\put(28,12.8){4}
\put(19,21.8){5}
\put(11,7.5){\line(0,1){4}}
\put(29,7.5){\line(0,1){4}}
\put(12.77,15.77){\line(1,1){5.46}}
\put(27.23,15.27){\line(-1,1){5.46}}
\end{picture}
\end{center}}
\end{minipage}
\caption{A decomposable submodels in a 5-way case}
\label{ex_fgg}
\end{figure}
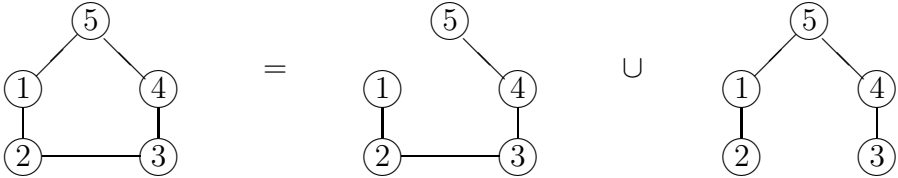

In this section, we compare our Algorithm 1 with the conventional IPS
by numerical experiments. %
We consider $J$-way cycle model with the generating class 
$\{ \{1,2\},\{2,3\},\dots, \{J-1,J\}, \{J,1\}\}$ for $J \ge 4$.
As a family of decomposable submodels which %
span the 
model we use the set of two decomposable submodels
obtained by deleting one element of generating class of the
hierarchical model.
We show the considered model and its submodels in Table \ref{submodels},
where $\{1,2\}$ is abbreviated as $12$. 
For example in the 5-way case we span
$M_5=\{12,23,34,45,15\}$ by $M_5\setminus \{15\}$ and $M_5\setminus \{23\}$
as illustrated in  Figure \ref{ex_fgg}.

Before we present the results of the experiments, 
we consider the space-saving implementation of Algorithm 1.

\begin{table}[htb]
\caption{The submodels in numerical experiments}
\begin{center}
\begin{tabular}{|c|c|c|}
\hline
Dim&Hierarchical model&Decomposable submodels\\
\hline
$4$ & $M_4=\{12,23,34,14\}$ & $M_4\setminus \{14\},M_4\setminus \{23\}$\\
\hline
$5$ & $M_5=\{12,23,34,45,15\}$ & $M_5\setminus \{15\},M_5\setminus \{23\}$\\
\hline
$6$ & $M_6=\{12,23,34,45,56,16\}$ & $M_6\setminus \{16\},M_6\setminus \{34\}$\\
\hline
$7$ & $M_7=\{12,23,34,45,56,67,17\}$ & $M_7\setminus \{17\},M_7\setminus \{34\}$\\
\hline
$8$ & $M_8=\{12,23,34,45,56,67,78,18\}$ & $M_8\setminus \{18\},M_8\setminus \{45\}$\\
\hline
\end{tabular}
\label{submodels}
\end{center}
\end{table}

\subsection{Tree-computation of Algorithm 1}
For the conventional IPS, the implementation of the tree-computation  
has been considered in Jirou\v{s}ek \cite{i11} and Jirou\v{s}ek and
P\v{r}eu\v{c}il \cite{i10}. 
Badsberg and Malvestuto\cite{badsberg} improved the algorithm by
applying the Markovian information propagation techniques with junction
trees of the triangulated models.

\begin{figure}[htpb]
 \centering
 \label{clique-tree}
 \includegraphics{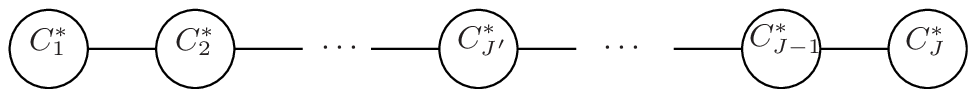}
 \caption{The junction tree for the $J$-way cycle model}
 \centering
 \label{directed-tree}
 \includegraphics{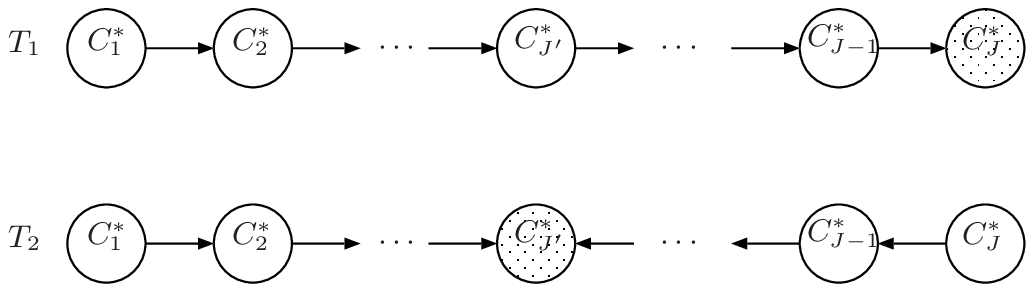}
 \caption{The directed trees}
\end{figure}

In this section we apply the Markovian propagation approach
to our Algorithm 1 for cycle models.
We triangulate the $J$-way cycle model by adding the edges 
$\{1,3\},\{1,4\},\ldots,\{1,J-1\}$.
Let ${\cal D}^*$ and ${\cal S}^*$ denote
the triangulated model and the set of %
separators of ${\cal D}^*$, 
$$
{\cal D}^* = \{\{1,2,3\},\{1,3,4\},\ldots,\{1,J-1,J\}\},\quad
{\cal S}^* = \{\{1,3\},\{1,4\},\ldots,\{1,J-1\}\}. 
$$
Then we note that the cell probabilities $p(i)$ satisfy 
$$
p(i) = \frac{\prod_{C \in {\cal D}^*} p(i_C)}
{\prod_{S \in {\cal S}^*} p(i_S)}.
$$
So we consider the tree-computation algorithm which updates
$q^{(t)}(i_C)$,  
$C \in {\cal D}^*$ instead of $q^{(t)}(i)$. 
Let $|{\cal I}_C^*| = \max_{C \in {\cal C}^*} |{\cal I}_C|$. 
While the computational cost per an update procedure of Algorithm 1 
is $O(|{\cal I}|)$, 
that of the tree-computation algorithm is reduced to 
$O(|{\cal I}_C^*|)$.
Denote by $C_j^*$ a generator %
$\{1,j-1,j\}$, $j=1,\ldots,J$. 
The junction tree for ${\cal D}^*$ is uniquely defined as in Figure 
\ref{clique-tree}. 
The decomposable submodels we use are 
$M^1_J = M_J \setminus \{J-1,J\}$ and 
$M^2_J = M_J \setminus \{J'-1,J'\}$ for some $1 < J' < J$. 
Direct the junction tree in two ways such that 
$C_J$ and $C_{J'}$ are the unique sink as in Figure \ref{directed-tree}
and denote them by $T_1$ and $T_2$, respectively. 
Then the Markovian propagation algorithm proposed here is described as 
information propagation on $T_1$ and $T_2$.

\ \\
\noindent {\bf Algorithm 4}\\
\quad Define 
$r^{(t)}(i_{\{1,2\}}) = r(i_{\{1,2\}})$ and 
$r^{(t)}(i_{\{J-1,J\}}) = r(i_{\{J-1,J\}})$ for all $t$\\
\ \\
\quad (1) Update via $M_J^1$ \\
\ \\
\quad \quad for $j = 3$ to $J$ do\\
\quad \quad \quad if $j \neq J$, update $q^{(t)}(i_{C_j^*})$ and
$r^{(t)}(i_{\{1,j\}})$ 
by 
\begin{equation}
 \label{update:q-1}
  q^{(t+1)}(i_{C_j^*}) = 
  q^{(t)}(i_{C_j^*}) \times 
  \frac{r^{(t+1)}(i_{\{1,j-1\}}) r(i_{\{j-1,j\}})}{r(i_{\{j-1\}})} 
  \times
  \frac{q^{(t)}(i_{\{j-1\}})}{q^{(t)}(i_{\{1,j-1\}}) \times
  q^{(t)}(i_{\{j-1,j\}})}, 
\end{equation}
$$
r^{(t+1)}(i_{\{1,j\}}) = \sum_{i_{j-1}} q^{(t+1)}(i_{C_j^*})
$$
\quad \quad \quad and send $r^{(t)}(i_{\{1,j\}})$ to $C_{j+1}^*$.\\
\quad \quad \quad if $j=J$, update $q^{(t)}(i_{C_J^*})$ by
(\ref{update:q-1}).\\
\ \\
\quad (2) Update via $M_{J}^2$ \\
\ \\
\quad \quad for $j = 3$ to $J'-1$ \\
\quad \quad \quad update $q^{(t)}(i_{C_j^*})$ and $r^{(t)}(i_{\{1,j\}})$
by (\ref{update:q-1})
and send $r^{(t)}(i_{\{1,j\}})$ to $C_{j+1}^*$.\\
\ \\
\quad \quad for $j = J$ to $J'+1$ \\
\quad \quad \quad update $q^{(t)}(i_{C_j^*})$ and 
$r^{(t)}(i_{\{1,j-1\}})$ by 
$$
q^{(t+1)}(i_{C_j^*}) = 
q^{(t)}(i_{C_j^*}) \times
\frac{r^{(t+1)}(i_{\{1,j\}}) r(i_{\{j-1,j\}})}{r(i_{\{j\}})} 
\times 
\frac{q^{(t)}(i_{\{j\}})}{q^{(t)}(i_{\{1,j\}}) 
q^{(t)}(i_{\{j-1,j\}})}, 
$$
$$
r^{(t+1)}(i_{\{1,j-1\}}) = \sum_{i_{j}} q^{(t+1)}(i_{C_j^*})
$$
\quad \quad \quad and send 
$r^{(t)}(i_{\{1,j-1\}})$
to $C_{j-1}^*$.\\  
\ \\
\quad \quad Update $q^{(t)}(i_{C_{J'}^*})$ by 
$$
q^{(t+1)}(i_{C_{J'}^*}) = 
q^{(t)}(i_{C_{J'}^*}) \times
\frac{r^{(t+1)}(i_{\{1,J'-1\}}) r^{(t+1)}(i_{\{1,J'\}})}
{r(i_{\{1\}})} 
\times
\frac{q^{(t)}(i_{\{1\}})}{q^{(t)}(i_{\{1,J'-1\}}) 
q^{(t)}(i_{\{1,J\}})}.
$$
\hfill\qed\\

It is easy to show that 
\begin{align}
 \label{eq:update}
q^{(t+1)}(i) &= 
\frac{\prod_{C \in {\cal D}^*} q^{(t+1)}(i_C)}
{\prod_{S \in {\cal S}^*} q^{(t+1)}(i_S)} \notag\\ 
&= 
\frac{\prod_{C \in {\cal D}^*} q^{(t)}(i_C)}
{\prod_{S \in {\cal S}^*} q^{(t)}(i_S)} \times
\frac{\prod_{C \in {\cal C}_1} r(i_C)}
{\prod_{S \in {\cal S}_1} r(i_S)} \times
\frac{\prod_{S \in {\cal S}_1} q^{(t)}(i_S)} 
{\prod_{C \in {\cal C}_1} q^{(t)}(i_C)}\notag\\
&= 
q^{(t)}(i) \times
\frac{\prod_{C \in {\cal C}_1} r(i_C)}
{\prod_{S \in {\cal S}_1} r(i_S)} \times
\frac{\prod_{S \in {\cal S}_1} q^{(t)}(i_S)} 
{\prod_{C \in {\cal C}_1} q^{(t)}(i_C)}. 
\end{align}
(\ref{eq:update}) looks the same as (\ref{eq:a1}). 
However there are some minor differences.
In (\ref{eq:a1})  
$q^{(t)}(i_S) = \sum_{i_{\Delta \setminus S}} q^{(t)}(i)$.  
However $q^{(t)}(i_S)$ in (\ref{eq:update}) is derived by 
$q^{(t)}(i_S) = \sum_{i_{C \setminus S}} q^{(t)}(i_C)$ for some 
$C \in {\cal C}^*$. 
Since $q^{(t)}(i_S)$ is not necessarily normalized,
$\sum_{i_{C \setminus S}} q^{(t)}(i_C) \neq 
\sum_{i_{C' \setminus S}} q^{(t)}(i_{C'})$ 
for $C \neq C'$ and $C$, $C' \in {\cal C}^*$. 
Hence in general 
$\sum_{i_{C \setminus S}} q^{(t)}(i_C) \neq 
\sum_{i_{\Delta \setminus S}} q^{(t)}(i)$.
In this sense Algorithm 4 is an approximate algorithm for Algorithm 1.  

In the experiments, we compare the performance of Algorithm 4 and 
the Markovian propagation
algorithm for the conventional IPS by Badsberg and Malvestuto\cite{badsberg}. 
\subsection{The results of the numerical experiments}
In this section we present the results of numerical experiments.
We set $I_1 = \cdots =I_{\delta} =I$ and $I=2$, $3$ or $4$.
We generated random contingency tables by filling each cell by uniform
random integers from 1 to $10^6$ 
and we obtained MLE by Algorithm 4 and the Markovian propagation
algorithm for the conventional IPS by Badsberg and
Malvestuto\cite{badsberg}. 
As the convergence criterion we used 
$$
\sum_{j=3}^J \sum_{i_{C_j^*} \in {\cal I}_{C_j^*}}
|q^{(t+1)}(i_{C_j}^*)-q^{(t)}(i_{C_j}^*)|
\leq 10^{-6}.  
$$
For each dimension and each number of levels, 
we generated 1000 contingency tables and took 
the average of the CPU time and the number of steps to convergence.    
Denote by $\tau$ and $\tau_{\text{conv}}$ the CPU time for Algorithm 4
and the conventional IPS, respectively.
Let $\nu$ and $\nu_{\text{conv}}$ be the number of steps to
convergence for Algorithm 4 and the conventional IPS, respectively.
We also calculated the probability that $\tau < \tau_{\text{conv}}$ and 
$\nu/\nu_{\text{conv}}$. 
The computation was done on a Pentium IV 3.2GHz CPU machine.

The results are shown in Table \ref{pic3}.  
In all of our runs Algorithm 4 converged to MLE.  
The experiments show that Algorithm 4 converges faster when the
dimension is larger than 7. 
The computational cost per an update of Algorithm 4 is expected to be
larger than that of the conventional IPS.
As we can see from Table \ref{pic3}, however, the number of steps to
convergence of Algorithm 4 is smaller than that of the conventional IPS.
$\nu/\nu_{\text{conv}}$ gets smaller as the dimension of the model 
gets larger.
Therefore the results of the experiments suggest that Algorithm 4 is
more efficient than the conventional IPS when the dimension of the model
is large for general hierarchical models.

\begin{table}[htb]
 \centering
 \caption{CPU time and the number of steps to convergence}
 \label{pic3}
 (i) $I=2$\\
 \begin{tabular}{|c|cc|c|cc|c|}\hline
  & \multicolumn{2}{|c|}{CPU time} & 
  $\mathrm{Pr}(\tau < \tau_{\text{conv}})$
  & \multicolumn{3}{|c|}{number of steps}\\ \cline{2-7}
  Dim & $\tau_{\text{conv}}$ & $\tau$ & & $\nu_{\text{conv}}$ & $\nu$ & 
  $\nu/\nu_{\text{conv}}$\\ \hline
  4   & 0.0156 & 0.0178 & 0.171 & 11.652 & 6.887 & 0.591\\
  5   & 0.0193 & 0.0205 & 0.237 &  9.391 & 4.413 & 0.470\\
  6   & 0.0233 & 0.0198 & 0.465 &  7.841 & 3.348 & 0.427\\
  7   & 0.0289 & 0.0204 & 0.884 &  8.000 & 3.000 & 0.375\\
  8   & 0.0407 & 0.0258 & 0.957 &  9.000 & 3.000 & 0.333\\ \hline
 \end{tabular}\\
\ \\
 (ii) $I=3$\\
 \begin{tabular}{|c|cc|c|cc|c|}\hline
  & \multicolumn{2}{|c|}{CPU time} & 
  $\mathrm{Pr}(\tau < \tau_{\text{conv}})$
  & \multicolumn{3}{|c|}{number of steps}\\ \cline{2-7}
  Dim & $\tau_{\text{conv}}$ & $\tau$ & & $\nu_{\text{conv}}$ & $\nu$ & 
  $\nu/\nu_{\text{conv}}$\\ \hline
  4   & 0.0337 & 0.0440 & 0.023 & 11.098 & 6.428 & 0.579\\ 
  5   & 0.0455 & 0.0463 & 0.394 &  9.345 & 4.455 & 0.477\\
  6   & 0.0484 & 0.0451 & 0.469 &  7.000 & 3.000 & 0.429\\
  7   & 0.0697 & 0.0559 & 0.929 &  8.000 & 3.000 & 0.375\\
  8   & 0.0951 & 0.0672 & 0.997 &  9.000 & 3.000 & 0.333\\ \hline
 \end{tabular}\\
\ \\
 (iii) $I=4$\\
 \begin{tabular}{|c|cc|c|cc|c|}\hline
  & \multicolumn{2}{|c|}{CPU time} & 
  $\mathrm{Pr}(\tau < \tau_{\text{conv}})$
  & \multicolumn{3}{|c|}{number of steps} \\ \cline{2-7}
  Dim & $\tau_{\text{conv}}$ & $\tau$ & & $\nu_{\text{conv}}$ & $\nu$ & 
  $\nu/\nu_{\text{conv}}$\\ \hline
  4   & 0.0665 & 0.0941 & 0.000 & 10.493 & 4.943&  0.471\\ 
  5   & 0.0722 & 0.1032 & 0.041 &  7.080 & 2.980&  0.421\\
  6   & 0.1005 & 0.1007 & 0.324 &  7.000 & 3.000&  0.429\\
  7   & 0.1437 & 0.1254 & 0.971 &  8.000 & 3.000&  0.375\\
  8   & 0.2028 & 0.1551 & 0.997 &  9.000 & 3.000&  0.333\\ \hline
 \end{tabular}
\end{table}

\section{Some discussions}
\label{sec:discussions}

For using the proposed algorithms, we have to find a family of
decomposable submodels that %
span a generating class of a
hierarchical model.  We recommend spanning the generating class by
a small number of large decomposable submodels.  Here large  decomposable
submodels might mean maximal submodels in the sense of model inclusion
or submodels with largest degrees of freedom.
In the literature some methods for finding a
maximal chordal subgraph of a given graph are studied (\cite{g4},
\cite{g6}, \cite{g5}). In the case of graphical models, this might
give a solution to our problem. However we have
to satisfy the condition that each element of a generating class 
is contained in at least one decomposable submodel.  Therefore we 
need a method to find a maximal chordal subgraph under the
restriction that specific generators %
are contained.

A referee suggested the following simple algorithm.
Suppose that a model $M({\cal C})$ with $|{\cal C}| = m$ is given.
For each set $C$ in ${\cal C}$\\
\hspace*{1cm} choose an ordering $(C_1, \dots  C_m)$
of sets in $\cal C$ such that $C_1 = C$ and 
$|C_{j-1} \cap  C_j| = \max$\\ \hspace*{1.5cm} for all $j \ge 2$; \\
\hspace*{1cm} set $ C' = \{C_1\}$;\\
\hspace*{1cm} for $j = 2,\dots,m$\\
\hspace*{2cm} if ${\cal C}' \cup \{C_j\}$ has the running intersection
property then set ${\cal C}' := {\cal C}' \cup \{C_j\}$. \\
Note that testing the running intersection property on a set family
takes linear time \cite{ty}.

In this paper we compared various algorithms of IPS in terms of the
CPU time to convergence. 
We showed that proposed algorithm converges faster than conventional IPS
when the model is large by numerical experiments.
We consider the implementation of the tree-computation of Algorithm 1
only in the case of cycle models. 
It may be possible to implement the tree-computation of Algorithm 1
for general hierarchical model when the decomposable submodels are
given. 
This topic needs further investigation and is left to our future research.

\bigskip
\noindent
{\bf Acknowledgment.} \quad  The authors are grateful to Hisayuki Hara
for implementation of the tree-computation in Section 4 and to 
Satoshi Kuriki for very useful comments. They thank two referees for
very constructive and detailed comments.

\end{document}